\documentclass[conference]{IEEEtran}

\ifCLASSINFOpdf
\else
\fi

\usepackage{multirow}
\usepackage{latexsym}
\usepackage{amsmath}
\usepackage{amsfonts}
\usepackage{amssymb}
\usepackage{graphicx}
\usepackage{times}
\usepackage{cite}
\usepackage{kbordermatrix}

\begin{document}

\title{Centralized Load Shedding Based on Thermal Limit of Transmission Lines Against Cascading Events}

\author{\IEEEauthorblockN{Bakhtyar Hoseinzadeh}
\IEEEauthorblockA{Dept. of Energy Technology, \\
Aalborg University, Denmark\\
Dept. of Electrical Engineering,\\
University of Kurdistan, Iran
}
\and
\IEEEauthorblockN{M. Hadi Amini}
\IEEEauthorblockA{Dept. of Elec. and Com. Eng.,\\
Carnegie Mellon University, Pittsburgh\\
PA 15213, USA\\
Yat-sen University- Carnegie melon University\\
Joint Institute of Engineering,\\
Guangzhou, Guangdong, China
}
\and
\IEEEauthorblockN{Claus Leth Bak}
\IEEEauthorblockA{Dept. of Energy Technology\\
Aalborg University, Denmark
}}

\maketitle

\begin{abstract}
Load shedding is the last and most expensive control action against system collapse and blackout. Achievement of an efficient emergency control to stabilize the power system following severe disturbances, requires two key objectives. First, preventing of further cascading outages, i.e. saving the available and determinant power system elements and second, issuing proper control actions to stabilize the power system inside the permissible time frame. In this paper online contingency analysis is performed to monitor secure and reliable operation of transmission lines. Load shedding locations are continuously updated based on loading rate/thermal limit of lines to prevent their outage. Simulation of severe contingencies carried out on 39 bus IEEE standard test system in DIgSILENT PowerFactory validates the efficiency of proposed method.
\end{abstract}

\IEEEpeerreviewmaketitle

\section{Introduction}
Blackouts rarely happen in the power system, but it does not diminish their significance. Power system design should be robust/resilient against uncertainties, unmodeled dynamics and should be able to withstand/tolerate wide range of combinational contingencies and cascading events to protect the stability of whole or as much as possible of power system. To this aim, proper preventive and self-healing control actions should be considered for different condition of power system.

In large scale power systems with thousands of complicated electrical equipments, although the current state variables may be available online for monitoring in the control center, the operators may not be able to make complicated decisions during serious contingencies, especially when the state variables are relatively fast comparing to the available time frame for bringing the variables back to the range above minimum permissible boundaries \cite{505}.

Load shedding is the last and most costly countermeasure/barrier against blackouts, which should be minimized according to the cost function of control action  \cite{508,103}. Load shedding is basically initiated/activated, if available primary and/or secondary control actions fail in stabilizing the system in their corresponding time span. Therefore, beside of tolerable range for key state variables of power system, e.g. frequency and voltage, time is also an important and determinant factor especially for protection system.

Most of recently proposed methods for centralized load shedding are merely based on frequency dynamics \cite{508,510,523}. The number of papers, in which voltage dynamics are also included are not considerable \cite{506,113,507,103,102}. Monitoring the frequency and voltage dynamics of power system is essential in the stability study, but they do not reflect reliable operation of determinant resources such as transmission lines. In other words, although a highly overloaded line may temporarily transmit the required active and reactive power causing the voltage and frequency profiles inside their permissible boundary, the overloaded line may not be available for a long time beyond the tolerable time of their protection relay. If the loading rate of a highly overloaded line is not alleviated in time, their over current protection relay may be triggered resulting in outage of line. Under such a circumstances even though the state variables seem to be in the safe range, they suddenly and unexpectedly collapse leading to initiate the next successive cascading events \cite{103,102}.

The current amplitude of transmission lines are available in the control center using existing Wide Area Measurement System (WAMS) collecting data from Phasor Measurement Unit (PMU) or Synchrophasors typically installed at substations in the transmission level.
\section{Monitoring of Loading Rate of Transmission Lines}\label{sec:3}
The thermal overload protection provides tripping
or alarming based on a thermal model of line currents. The curve and time dial settings of inverse time-overcurrent characteristic of protection relays conform to the \textit{IEEE C37.112} and \textit{IEC 255-03} standards \cite{1023}. The maximum available time before activation of over current relay of line $k$ can be calculated based on the rate of thermal overload and the following inverse time tripping characteristic:
\begin{equation}\label{eq:1}
T_k^{oc}(t) = \frac{\alpha }{{{({r_k^{oc}(t)})}^{\gamma}} - 1} + \beta
\end{equation}
Parameters $\alpha$, $\beta$ and $\gamma$ are constants to provide desired curve characteristics from aforementioned standards according to the application. $r_k^{oc}(t)$ represents the over current rate regarding the chosen current setting/pickup ($i_k^{s}$) of the relay, which should be greater/less than one for trip/reset operation of the relay:
\begin{equation}\label{eq:2}
r_k^{oc}(t)=\frac{i_k(t)}{i_k^{s}} 
\end{equation}

\begin{figure}[b]
	\centering
	\includegraphics[width=0.7\linewidth]{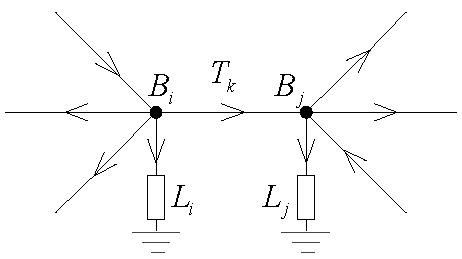}
	\caption{Typical bus-line topology}
	\label{fig:1}
\end{figure}
For mathematical formulation of proposed method, a network with $b$ Transmission lines ($T_1$-$T_b$) and $n$ Buses ($B_1$-$B_n$) is considered.
Fig.~\ref{fig:1} indicates a typical bus-line connection with power flow direction regardless of their existing electrical elements. $B_i$, $L_i$ and $T_k$ stand for bus $i$, Load $i$ and line $k$, respectively.

Outage risk index of an overloaded line is defined based on its loading rate to prevent more cascading events. To decline the line overload, e.g. $T_k$ in Fig.~\ref{fig:1}, only alleviation of output power flow of destination bus ($B_j$) may be helpful, i.e the transmission line/s with outbound power flow and existing load/s ($L_j$). Therefore, the load/s of destination bus ($B_j$) should be assigned higher priority for interruption. To this aim, power flow direction and loading rate ($i_k^{oc}(t)$) are considered.

The power flow direction of lines, e.g. $T_k$, can be summarized in the incident matrix ($A(t)$), which describes the network topology, i.e. the connection between buses and lines:
\begin{equation}\label{eq:3}
A_{n \times b}(t) = \kbordermatrix{
         	&       1& \hdots  & k         & \hdots  & b\\
	\vdots  & \hdots & \hdots  & 0         & \hdots  & \hdots\\
       	i   & \hdots & \hdots  & +{D_k}(t) & \hdots  & \hdots\\
	\vdots  & \hdots & \hdots  & 0         & \hdots  & \hdots\\
	     j  & \hdots & \hdots  & -{D_k}(t) & \hdots  & \hdots\\
	\vdots  & \hdots & \hdots  & 0         & \hdots  & \hdots
}
\end{equation}
where ${D_k}(t)$ represents the power flow Direction (${P_k}(t)$) of line $T_k$:
\begin{equation}\label{eq:4}
A_{i,k}(t)=D_k(t) \buildrel \Delta \over =  \left\{ {\begin{array}{*{20}{c}}
	+1&{{P_k}(t)\text{ enters bus } i}\\
	{-1}&{{P_k}(t)\text{ leaves bus } i}\\
	~~0&\text{ No incident with bus } i
	\end{array}} \right.
\end{equation}
The vector $i^{oc}(t)$ represents the loading rate of lines defined in (\ref{eq:2}):
\begin{equation}\label{eq:5}
i^{oc}(t)=\left[ {\begin{array}{*{20}{c}}
	{i_1^{oc}(t)}\\
	\vdots \\
	{i_k^{oc}(t)}\\
	\vdots \\
	{i_b^{oc}(t)}
	\end{array}} \right]
\end{equation}
The Impact Factor (IF) of the load installed at bus $i$ on its transmission lines considering both power flow direction ($A_{i,k}(t)$) and line loading rate ($i_k^{oc}(t)$) can be determined as follows:
\begin{equation}\label{eq:6}
IF_{i,k}(t)=A_{i,k}(t) \cdot i_k^{oc}(t)
\end{equation}
Assume that $c$ is the index of the line with the maximum Impact Factor (IF) from the load $i$ over the range of $k\in[1,...,b]$:
\begin{equation}\label{eq:10}
IF_{i,c}(t) = {\mathop {max}\limits_{k = 1}^b (IF_{i,k}(t))}
\end{equation}
To sort the network loads and therefore to prioritize their interruption order, the maximum available time before outage of their overloaded transmission line/s resulting from aforementioned loading rate is considered. To this aim, the inverse trip time calculated in (\ref{eq:1}) associated with the most overloaded line, i.e. the line with the index $c$ is assigned to the corresponding load bus:
\begin{equation}\label{eq:11}
T_i^{oc}(t) = T_c^{oc}(t)
\end{equation}

\begin{figure}[b]
	\centering
	\includegraphics[width=1\linewidth]{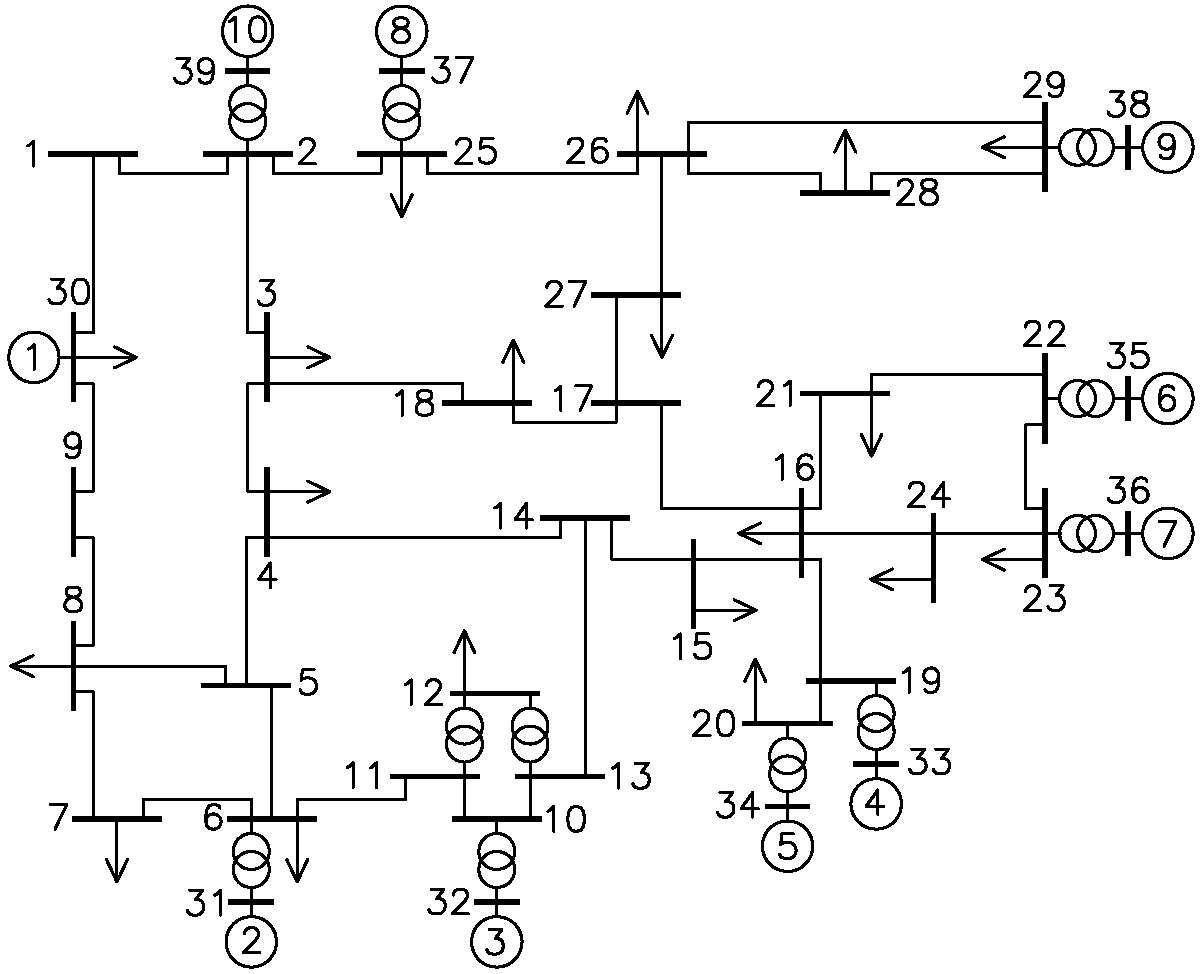}
	\caption{39 bus IEEE standard test system}
	\label{fig:4}
\end{figure}

\section{Simulation Setup}\label{sec:11}
In order to evaluate the performance of proposed scheme, the 39 bus IEEE standard test system demonstrated in Fig.~\ref{fig:4} is chosen as case study in DIgSILENT PowerFactory software \cite{106,512,107}. Synchronous generators are equipped with IEEE standard governor GOV-IEESGO and automatic voltage regulator AVR-IEEEX1. Moreover, dependency of active and reactive power of loads to voltage and frequency is considered \cite{a1,a2,a3,ashk}. The loading rate of transmission lines in percentage of their nominal current are assumed to be equal to 85 \%.

\section{Simulation Results}\label{sec:12}

\begin{figure}[!b]
	\centering
	\includegraphics[width=1.0\linewidth]{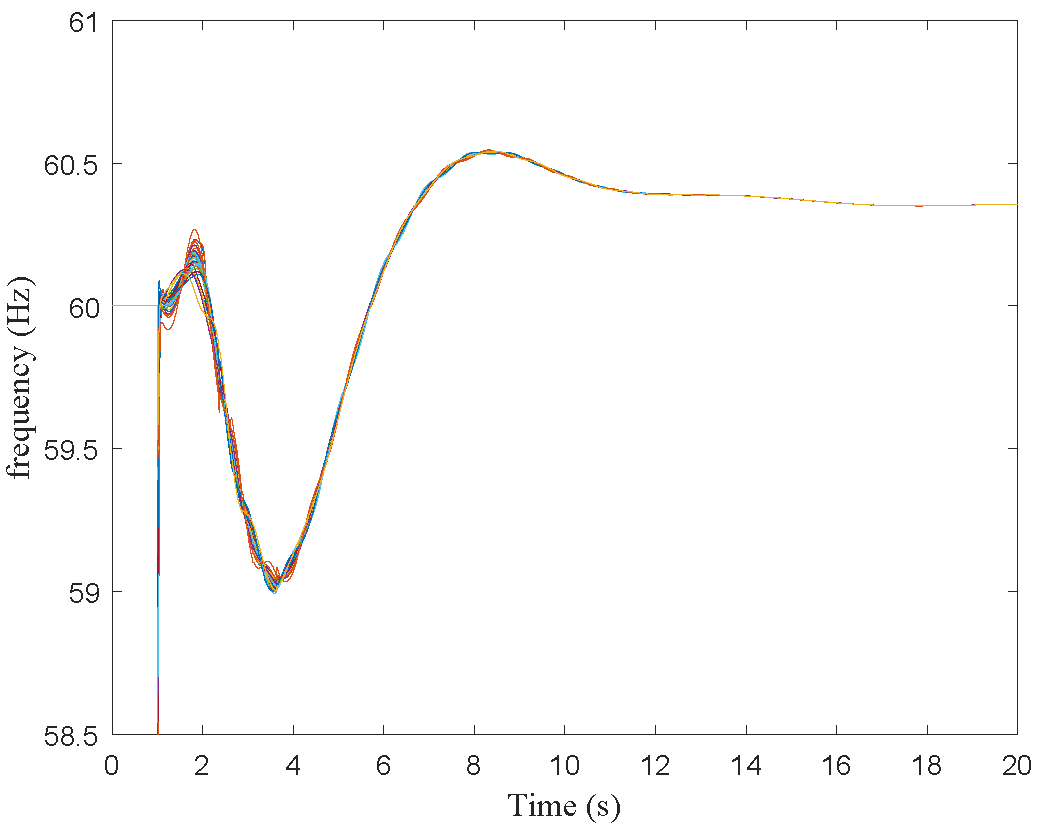}
	\caption{Frequency at different buses}
	\label{fig:5}
\end{figure}

\begin{figure}[!b]
\centering
\includegraphics[width=1.0\linewidth]{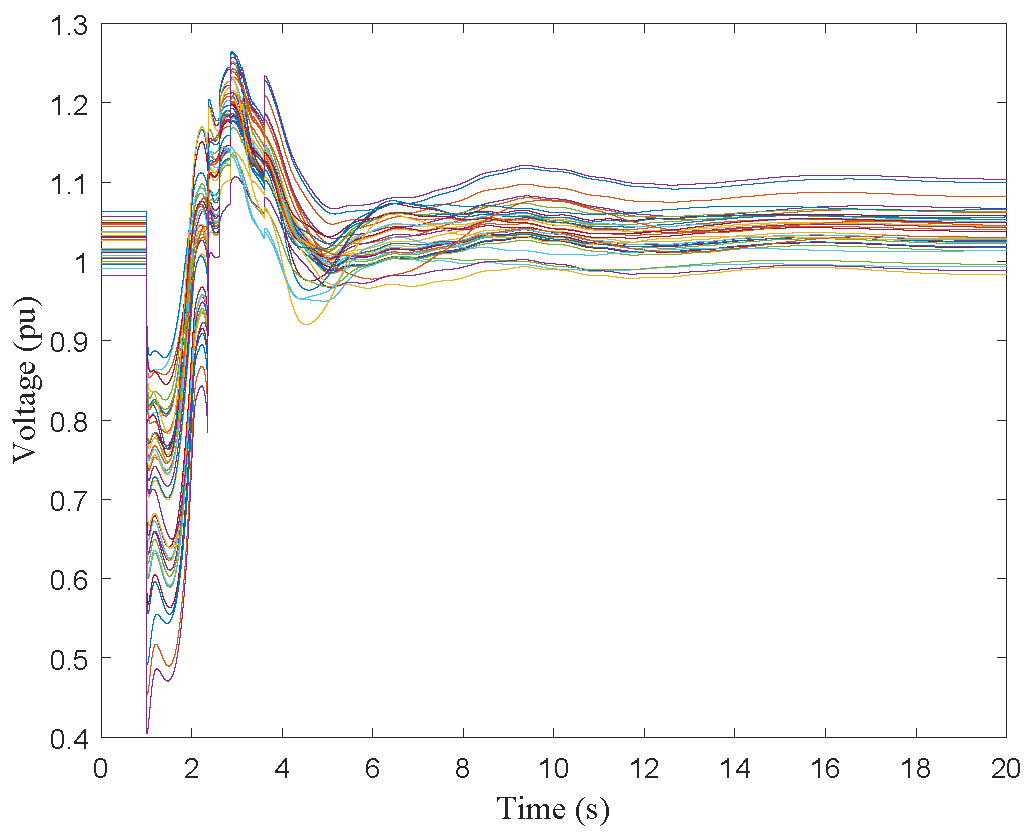}
\caption{Voltage at different buses}
\label{fig:6}	
\end{figure}

\begin{figure}[!b]
\centering
\includegraphics[width=1.0\linewidth]{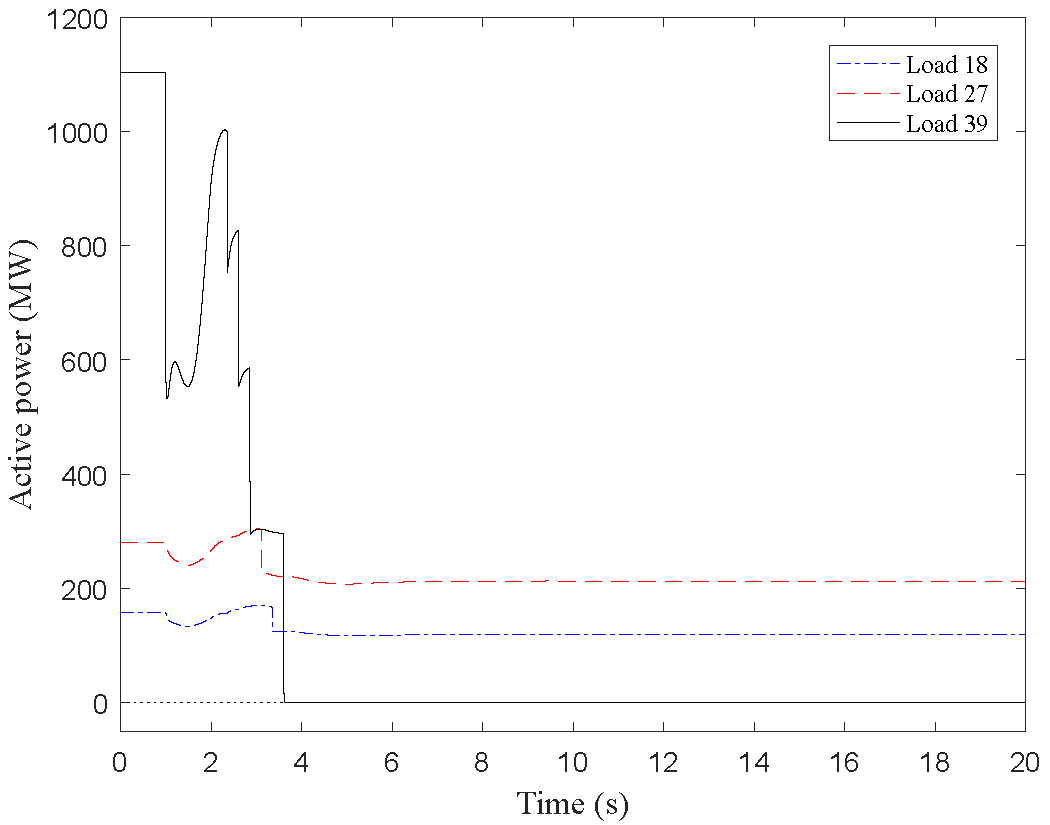}
\caption{Affected Loads by load shedding scheme}
\label{fig:7}	
\end{figure}

\begin{figure}[!b]
\centering
\includegraphics[width=1.0\linewidth]{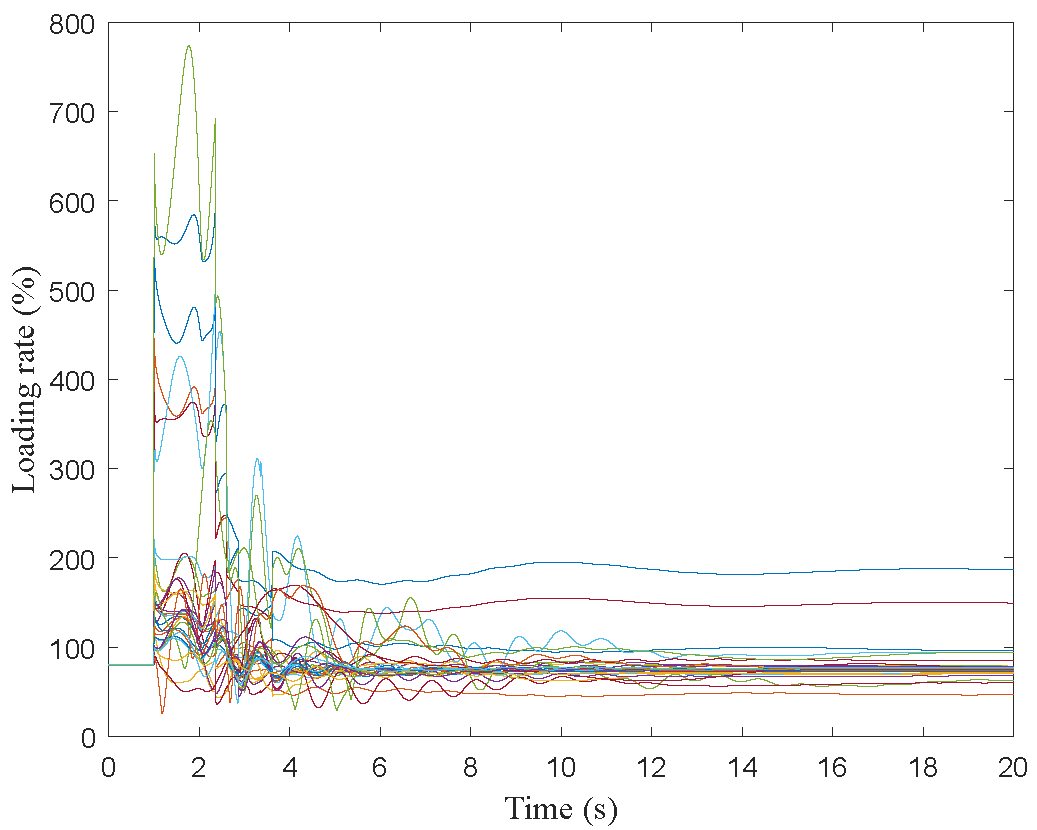}
\caption{Loading rate of transmission lines}
\label{fig:8}
\end{figure}

Outage of generator 1 at second 1 connected to the bus 30 with the generation almost equal to 1012 MW (17.7 \% of total generation) is chosen as a severe contingency (Fig.~\ref{fig:4}). Figs.~\ref{fig:5}-\ref{fig:8} indicates the simulation results for 20 s. Before the disturbance, the power system is in steady state with the frequency stable at 60 Hz (Fig.~\ref{fig:5}) and bus voltages in range around 1 pu of nominal (Fig.~\ref{fig:6}).
Figs.~\ref{fig:7} and \ref{fig:8} show the active power of loads and loading rate of transmission lines (equal to 85\% prior the event), respectively.

Outage of the largest generator of the system, i.e. G1 at 1s causes a sharp plunge at frequency and voltage profiles. Since the frequency is generally smoother than voltage dynamics, due to the existing inertia of the system, frequency returns back to a value close to its nominal value. Meanwhile, some minor local or inter-area oscillation modes are also observable in frequency curve \cite{521,109}.

Sudden voltage drop at nearby buses from 1 to 0.4 pu (Fig.~\ref{fig:6}) causes reduction of load's active power, e.g. load 30 form 1100 to 550 MW as can be seen in Fig.~\ref{fig:7}. The dependency of load's active and reactive power to the voltage dynamics justifies this phenomenon, which prevents fast drop of frequency as it is observed at voltage profile. Depending on the composition of loads, it may even improve the frequency temporarily as it can be seen around 2 s.

Reduction of load's active power consequently alleviates the load-generation imbalance, which causes relatively recovery of bus voltages. Recovery of voltages and thereafter load's active power brings the load-generation imbalance back to its real situation. It means that frequency starts to fall down, even though the voltages are equal or above nominal values.

Under such a condition, the system suffers form a frequency event rather than a voltage stability issue. The demanded active power is transfered to the event location form far away generation units, which increases the loading burden of transmission lines on the way to the destination. Fig.~\ref{fig:8} indicates that the loading rate of lines connected to the buses 1, 8 and 9 are increased to a high value between 6 to 8 times of their nominal current, which may trigger their over current protection relay and hence sudden outage of them.

The proposed technique continuously monitors the transmission line currents to identifies the loads causing overloading of lines. The loads connected to the end of most overloaded lines are assigned highest priorities for curtailment. As can be seen the loading rate of line pron to outage are reduced by proper load shedding decision in time at the best locations before operation of corresponding protection relay.

According to Fig.~\ref{fig:7}, the first three feeders of load 30 are successively disconnected by three stages of load shedding relay followed by interruption of the first stage of loads 27, 18 and the last stage of load 30 to bring the frequencies, voltages and loading rates back to the acceptable range.

\section{Conclusion}
To prevent cascading outage of transmission lines following contingencies due to operation of their over current protection relays, their thermal limit criterion is considered in the load shedding control action. Proper load shedding locations are determined online using performing the contingency analysis based on thermal limit of transmission lines. The loads supplied by most overloaded transmission lines are shed first.

\end{document}